\documentclass[10pt]{amsart}
\usepackage{amssymb}

\theoremstyle{plain}
\newtheorem{theorem}{Theorem}[section]

\newtheorem{lemma}[theorem]{Lemma}
\newtheorem{cor}[theorem]{Corollary}

\newtheorem{rem}[theorem]{Remark}

\newcommand{\Z}{{\mathbb Z}}
\newcommand{\Q}{{\mathbb Q}}
\newcommand{\C}{{\mathbb C}}
\newcommand{\R}{{\mathbb R}}
\newcommand{\K}{{\mathbb K}}

\numberwithin{equation}{section}

\begin{document}

\title[Common values]{Common values of a class of linear recurrences}

\author[Attila~Peth\H{o}]{Attila~Peth\H{o}}
\address{Department of Computer Science,\\ University of Debrecen,\\
H-4002 Debrecen, P.O. Box 400}
\email{petho.attila@unideb.hu}

\email{}

\begin{abstract}
Let $(a_n), (b_n)$ be linear recursive sequences of integers with characteristic polynomials $A(X),B(X)\in \Z[X]$ respectively. Assume that $A(X)$ has a dominating and simple real root $\alpha$, while $B(X)$ has a pair of conjugate complex dominating and simple roots $\beta,\bar{\beta}$. Assume further that $\alpha/ \beta$ and $\bar{\beta}/\beta$ are not roots of unity and $\delta = \log |\alpha|/ \log |\beta| \in \Q$. Then there are effectively computable constants $c_0,c_1>0$ such that the inequality
$$
  |a_n - b_m| > |a_n|^{1-(c_0 \log^2 n)/n}
$$
holds for all $n,m \in \Z^2_{\ge 0}$ with $\max\{n,m\}>c_1$.

\end{abstract}

\date{\today}

\maketitle

\section{Introduction}

Let $p_0,\ldots,p_{k-1}, g_0,\ldots,g_{k-1}\in \Z$ and
\begin{equation}\label{genrek}
  g_{n+k} = p_{k-1}g_{n+k-1}+\dots+p_0g_n
\end{equation}

for all $n\ge 0$ then $(g_n)$ is called {\it linear recursive sequence} with {\it characteristic polynomial} $G(X)=X^k-p_{k-1}X^{k-1}-\dots-p_0$. To avoid improper examples we assume $p_0$ and at least one of $g_0,\dots,g_{k-1}$ non-zero. Denote $\gamma_1,\dots,\gamma_l$ the distinct zeroes of $G(X)$ with multiplicities $m_1,\dots,m_l$ respectively. Set $\K=\Q(\gamma_1,\dots,\gamma_l)$. It is well known, see e.g \cite{ShT} Chapter C, that there exist polynomials $G_j(X)\in \K[X]$ of degree at most $m_j-1, j=1,\dots,l$ such that
$$
g_n = G_1(n)\gamma_1^n + \dots + G_l(n)\gamma_l^n
$$
holds for $n\ge 0$.

The breakthrough works of A. Baker \cite{B1}, \cite{B2} on lower bounds for linear forms in logarithms of algebraic numbers made it possible, among others, the systematic study of diophantine properties of linear recursive sequences. Moreover the results obtained with the help of Baker's theory are effective, they serve computable upper bounds for the size of the solutions. You find an excellent overview of earlier results in the book of Shorey and Tijdeman \cite{ShT}.

The subspace theorem of W.M. Schmidt \cite{Sch} became an other strong tool in the theory of linear recurrences, see e.g. J.-H. Evertse \cite{Evertse}. Usually the results obtained with the derivations of the subspace theorem are more general, but not effective.

One of the basic questions is the description of the common terms of recurrences. This includes values appearing several times in a recurrence. As recurrences are linear combinations of exponential functions, one expects that they have usually only finitely many common terms, and, in the opposite case, they have to be strongly related. This was proved ineffectively by J.-H. Evertse \cite{Evertse} and M. Laurent \cite{Laurent}. Effective version is known only for the equations $g_n = g_m$ provided the characteristic polynomial of $(g_n)$ has at most three roots with the same absolute value and for $g_n = h_m$ when both recurrences have dominating simple and real roots. The book of T.N. Shorey and R. Tijdeman \cite{ShT} gives detailed overview on the related results.

Recently E.F. Bravo et al. \cite{BGL}, \cite{BGLTK} established all solutions $n,m\in \Z$ of the diophantine equation $T_n = T_m,$ where $(T_n)$ denotes the tribonacci sequence, which is defined by the initial terms $T_{-1}=T_0=0, T_1=1$ and by the recursion $T_{n+3}=T_{n+2}+T_{n+1}+T_n, n\ge -1$. The relation \eqref{genrek} allows to extend $(g_n)$ to negative indices, but that terms are usually not integers. If, however, $p_0=\pm1$ the terms with negative indices are integers too. This holds to the tribonacci sequence, which negative branch is defined as $T_{-n-3}=-T_{-n-2}-T_{-n-1}+T_{-n}$. The characteristic polynomial of $(T_n)$ is $T(X) = X^3-X^2-X-1$, which has a dominating real root. On the other hand the characteristic polynomial of $(T_{-n})$ is $X^3+X^2+X-1=-X^3T(1/X)$, whose conjugate complex roots are dominating, i.e. the absolute values of them are larger than the absolute value of the third root. In this circumstances the older methods are not applicable, the authors needed a novel idea. As I mentioned Bravo et al. were interested to find all common terms for a given recurrence, but after a closer look became clear that the method is capable to prove effective finiteness result for a wide class of recurrences. My original goal was to generalize the result of E.F. Bravo et al. \cite{BGL}, \cite{BGLTK} to common values of $k$-generalized Fibonacci numbers, see e.g. \cite{P1}. Although I was unable to fulfill the original expectation, the results are general enough to summarize in these notes.

To formulate our theorems we have to introduce some notation. Let $(a_n), (b_n)$ be linear recursive sequences of integers with characteristic polynomials $A(X),B(X)\in \Z[X]$ respectively. Let $\alpha_1,\ldots,\alpha_k$ and $\beta_1,\ldots,\beta_l$ be the distinct roots of $A(X)$ as well as of $B(X)$ with multiplicities $m_1,\ldots,m_k$ and $n_1,\ldots,n_l$ respectively. We assume that $|\alpha_1|>|\alpha_2|\ge |\alpha_3|\ge\ldots\ge |\alpha_k|$; further $\beta_2=\bar{\beta}_1$ and  $|\beta_1|=|\beta_2|>|\beta_3|\ge\ldots \ge |\beta_l|$.

In the sequel $c_0,c_1,\ldots$ denote effectively computable constants depending only on the degrees and on the roots of the characteristic polynomials of $(a_n)$ and $(b_n)$. The algebraic numbers $\alpha, \beta$ are called multiplicatively dependent if there are non-zero integers $a,b$ such that $\alpha^a \beta^b =1$.

\begin{theorem} \label{t:main}
Let $(a_n)$ and $(b_m)$ be as above. Assume that $m_1=n_1=1, \alpha_1/ \beta_1$ and $\beta_2/\beta_1$ are not roots of unity and $\delta = \log |\alpha_1|/ \log |\beta_1| \in \Q$, i.e. $\alpha_1$ and $|\beta_1|$ are multiplicatively dependent. Then the diophantine inequality
\begin{equation}\label{e:alap}
  |a_n - b_m| > |a_n|^{1-(c_0 \log^2 n)/n}
\end{equation}
holds for all $n,m \in \Z^2_{\ge 0}$ with $\max\{n,m\}>c_1$, where
$$
c_0 > \frac{9c_{15}}{\log |\alpha_1|}
$$
and $c_1, c_{15}$ denote effectively computable positive constants defined in the proof.
\end{theorem}

An immediate consequence of our theorem is

\begin{cor} \label{c:main}
Under the assumptions of Theorem \ref{t:main} the equation 
$$
a_n = b_m
$$
has only finitely many, effectively computatble solutions in $n,m \in \Z^2_{\ge 0}$. 
\end{cor}

Using the theory of weighted $S$-unit equations, see e.g. \cite{Evertse_Gyory} it is easy to prove that if $\alpha_1/ \beta_1$ and $\beta_2/\beta_1$ are not roots of unity then the equation $a_n - b_m = K$ has for any given integer $K$ only finitely many solutions in $(n.m)\in \Z_{\ge 0}^2$. We do not give here the proof, but refer to the proof of Theorem 3.4 of \cite{P1}.

The assumption $\alpha_1/ \beta_1$ and $\beta_2/\beta_1$ are not roots of unity is natural because, otherwise, $a_n = b_m$ may have infinitely many solutions. Indeed, set $\alpha_1=2, \beta_1=2i, \beta_2=-2i$ and $a_n=\alpha_1^n, b_m=\beta_1^m + \beta_2^m$, then $b_{4k}=2^{4k+1}=a_{4k+1}$ hold for all $k\ge0$. Let now $(b_m)$ be the same and $(a_n)$ any sequence satisfying our assumptions and with $a_0=0$. As $b_{4k+1}=0$ for all $k\ge0$ the equation $a_n=b_m$ has infinitely many solutions.

On the other hand the assumption $\alpha_1$ and $|\beta_1|$ are multiplicatively dependent is not at al natural, moreover it is quite restrictive, but without it we can prove only finiteness, but usually not effective result. To highlight the reason consider the three appearances of Baker type inequalities in our proof, which are
\begin{enumerate}
  \item the proof of a lower bound on $|b_m|$,
  \item the proof of a lower bound on $|n-\delta m|$,
  \item the proof of an upper bound on $n$.
\end{enumerate}
The first two steps can be performed without any restrictions on $\delta$, but it seems to be essential for Step 3. 

Notice that Theorem \ref{t:main} requires properties only on the characteristic polynomials but not on the initial values of the sequences. Our second theorem gives a sufficient condition for pairs of linear recurrences for which Theorem \ref{t:main} is applicable. We show that the restriction $\log |\alpha_1|/\log |\beta_1|\in \Q$ narrow considerably the possible irreducible polynomials, but there are still infinitely many examples. Multiplying such polynomials with suitable ones we get the characteristic polynomials of linear recursive sequences for which the assumptions of Theorem \ref{t:main} hold. 
To avoid repetition of the phrase {\it maximal size of the roots of the polynomial $P(X)$} we denote it by $\overline{|P|}$.

\begin{theorem}\label{t:masodik}
  With $a,b,p,q\in \Z, p,q>0, \gcd(p,q)=1$ and $q$ even if $b$ is not a square define $Q_1(X)=X^2 + aX +b^p$ and $Q_2(X)=X^3+aX^2+bX +1$ such that $Q_2(X)$ has one real root outside the unit circle and a pair of conjugate complex roots. Assume that $P_1(X), P_2(X)\in \Z[X]$ and
\begin{itemize}
  \item[(i)] either\\
  $\bullet$ $a^2-4b^p<0,\;  \overline{|P_1|}<|b^{q/2}|, \;\overline{|P_2|}< \overline{|Q_1|}$,\\
  $\bullet$ $A(X)=(X \pm b^{q/2})P_1(X), \;B(X)=Q_1(X)P_2(X)$,
   \item[(ii)] or\\
  $\bullet$ $\overline{|P_1|}< \overline{|Q_2|}, \;\overline{|P_2|}< \overline{|X^3Q_2(1/X)|}$,\\
  $\bullet$  $A(X)=Q_2(X)P_1(X), \;B(X)=X^3Q_2(1/X)P_2(X)$.
\end{itemize}
Let $(a_n),(b_n)$ be non-zero linear recursive sequences with characteristic polynomials $A(X)$ and $B(X)$ respectively. If $\alpha$ and $\beta$ denote the dominating root of $A(X)$ and $B(X)$ respectively, then they satisfy the assumptions of Theorem 1.1. 
\end{theorem}

\begin{rem}
The assumption $Q_2(X)$ has one real root and a pair of conjugate complex roots means that the discriminant of $Q_2(X)$ is negative, i.e. $-27+18ab+a^2b^2-4a^3-4b^3<0$. 
\end{rem}

Combining Theorems 1.1 and 1.2 we obtain a far reaching generalization of Bravo et al. \cite{BGL}, \cite{BGLTK}. 

\begin{cor} \label{c:Bravo}
Let $a,b\in \Z$ with $-27-18ab+a^2b^2-4a^3+4b^3<0$ Let $f_0,f_1,f_2\in \Z$ not all zero and
$$
f_{n+3}=a f_{n+2}+b f_{n+1}+f_n,\; n\ge 0.
$$
Extend the sequence $(f_n)$ toward negative indices by the recursion
$$
f_{-n} = -b f_{-n+1} -a f_{-n+2} + f_{-n+3}, \; n\ge 1.
$$
Then there are only finitely many effectively computable $n,m\in \Z, n\not=m$ with $f_n = f_m$.
\end{cor}

\section{A Baker's type inequality}

An essential part of the proof of the main theorem depends on A. Baker's theory on linear forms in logarithms of algebraic numbers. Among the many variants we use here the one due to Matveev, which is one of the sharpest nowadays.

For an algebraic number $\eta$ with minimal polynomial
$$
f(X)=a_0(X-\eta^{(1)})\cdots(X-\eta^{(N)})\in {\mathbb Z}[X]
$$
with positive $a_0$, write $h(\eta)$ for its absolute logarithmic or Weil height given by
$$
h(\eta)=\frac{1}{N}\left(\log a_0+\sum_{j=1}^N \max\{0,\log |\eta^{(j)}|\}\right).
$$

We list important and well known properties of the function $h$ in the following lemma.

\begin{lemma}\label{l:Weil height}
Let $\gamma,\eta$ be algebraic numbers and $u\in \Q$. Then we have
\begin{enumerate}
  \item $h(\gamma \pm \eta) \le h(\gamma)+ h(\eta)+ \log 2$
  \item $h(\gamma \eta^{\pm1}) \le h(\gamma)+ h(\eta)$
  \item $h(\gamma^u) = |u| h(\gamma)$.
\end{enumerate}
\end{lemma}
For the proof see e.g. Waldschmidt \cite{W2000}, Ch. 3.2.

Let $\mathbb{K}$ be an algebraic number field of degree $d_{\mathbb{K}}$ and let $\eta_1, \eta_2, \ldots, \eta_t \in \mathbb{K}$  not $0$
or $1$, and $b_1, \ldots, b_t$ be nonzero integers. Put
$$
B =\max\{|b_1|, \ldots, |b_t|, 3\}\qquad{\rm and}\qquad
\Gamma = \prod_{i=1}^t \eta_i^{b_i}-1.
$$
Let $A_1, \ldots, A_t$ be positive integers such that
$$
A_j \geq \max \{d_{\mathbb{K}}h(\eta_j), |\log \eta_j|, 0.16\}, \quad {\text{\rm for}}\quad j=1,\ldots t.
$$
Under the circumstances above, Matveev \cite{M} proved
\begin{lemma}\label{Matveev}
 If $\Gamma \neq 0$, then
\begin{equation*}
\label{ineq:matveev} \log |\Gamma| > -3\cdot
30^{t+4}(t+1)^{5.5}d_{\mathbb{K}}^2(1+\log d_{\mathbb{K}})(1+\log tB)A_1A_2\cdots A_t.
\end{equation*}
\end{lemma}

\section{Proof of Theorem \ref{t:main}}

Under the assumptions there are polynomials $A_1,\ldots,A_k\in \Q(\alpha_1,\ldots,\alpha_k)[X]$, $B_1,\ldots,B_l \in \Q(\beta_1,\ldots,\beta_l)[X]$ of degrees $m_1-1,\ldots,m_k-1$ as well as $n_1-1,\ldots,n_l-1$ respectively such that
$$
a_n = A_1(n)\alpha_1^n+\ldots+A_k(n)\alpha_k^n, \quad b_n = B_1(n)\beta_1^n+\ldots+B_l(n)\beta_l^n
$$
hold for all $n\ge 0$. As $m_1=n_1=1$ and $\alpha_1$ and $\beta_1$ is dominating among the roots of $A(X), B(X)$ respectively we have $|\alpha_1|,|\beta_1|>1$, $\alpha_1$ is real $n_2=n_1=1$, $\beta_1 = \bar{\beta_1}$ and $A_1,B_1$ and $B_2$ are constants. Putting
$$
\delta_a= \frac{\log|\alpha_2|}{2\log |\alpha_1|} \quad \mbox{and} \quad \delta_b= \frac{\log|\beta_3|}{2\log|\beta_1|} 
$$
then $ \delta_a, \delta_b<1$ and it is an easy exercise to prove
\begin{eqnarray}
|a_n - A_1 \alpha_1^n| &<& c_2 |\alpha_1|^{\delta_a n}, \label{e:valoskorlat2}\\
c_3 |\alpha_1|^n &<& |a_n| < c_4 |\alpha_1|^n \label{e:valoskorlat3}
\end{eqnarray}
and
\begin{equation}\label{e:komplexkorlat2}
 |b_m - (B_1 \beta_1^m + B_2 \beta_2^m)| < c_5 |\beta_1|^{\delta_b m}
\end{equation}
with suitable positive effective constants $c_2,c_3,c_4,c_5$, provided $n$ and $m$ are large enough. 

In contrast, it is a consequence of the deep theory of A. Baker that if $\beta_1/\beta_2$ is not a root of unity then
\begin{equation}\label{e:komplexkorlat1}
 |\beta_1|^{m-c_6 \log m} < |B_1 \beta_1^m + B_2 \beta_2^m| < c_7 |\beta_1|^m,
\end{equation}
see Corollary 3.7 of \cite{ShT}. The last two inequalities imply
\begin{equation}\label{e:komplexkorlat3}
 |\beta_1|^{m-2c_6 \log m} < |b_m | < 2 c_7 |\beta_1|^m.
\end{equation}
provided $m$ is large enough. We derive only the lower bound, because the proof of the upper bound is straight forward
$$
|b_m | > |B_1 \beta_1^m + B_2 \beta_2^m| - c_8 |\beta_3|^m
  > |\beta_1|^{m-c_6 \log m} - c_8 |\beta_1|^{\delta_b m}
  > |\beta_1|^{m-2c_6 \log m}.
$$

\medskip

Contrary to the statement of Theorem \ref{t:main} we assume that there exist infinitely many $(n,m)\in \Z_{\ge 0}^2$ for which \eqref{e:alap} does not holds. We show that the choice of any such solutions with $n$ large enough leads to a contradiction. In the sequel let $(n,m)\in \Z_{\ge 0}^2$ be a solution of
\begin{equation*}
  ||a_n| - |b_m|| \le |a_n|^{1-(c_0 \log^2 n)/n}
\end{equation*}
and $n$ is large enough. As $||a_n| - |b_m|| = |a_n - b_m|,$ if $a_n$ and $b_m$ have the same sign and by $|-a_n - b_m|$, we have to deal with the two cases separately. In the sequel we are dealing with the first case because replacing $a_n$ by $-a_n$ in the argument below one gets easily the proof of the second case. Thus we are to do with the inequality
\begin{equation}\label{e:kisebb}
  |a_n - b_m| \le |a_n|^{1-(c_0 \log^2 n)/n}.
\end{equation}
We also assume $a_n\ge b_m$. The case $a_n < b_m$ can be handled similarly.

Instead of $a_n-b_m$, which may have many terms we will study  $A_{n,m}= A_1 \alpha_1^n - (B_1 \beta_1^m + B_2 \beta_2^m)$, which has only three summands, thus much more easy to treat. Moreover $|A_{n,m}|$ is not far from $a_n-b_m$. Indeed
$$
|A_{n,m} - (a_n - b_m)| \le |a_n - A_1 \alpha_1^n | + |b_m- (B_1 \beta_1^m + B_2 \beta_2^m)| < c_2 |\alpha_1|^{\delta_a n} + c_5 |\beta_1|^{\delta_a m}.
$$
Using \eqref{e:valoskorlat2} and \eqref{e:komplexkorlat2} we obtain
\begin{eqnarray}\label{e:kombinalt}
  |A_{n,m}| &\le& |a_n - b_m|+ |a_n - A_1 \alpha_1^n| + |b_m - (B_1 \beta_1^m + B_2 \beta_2^m)| \nonumber \\
  &\le&  |a_n|^{1-c_0(\log^2n)/n} + c_4 |\alpha_1|^{\delta_a n} + c_5 |\beta_1|^{\delta_b m} \\
  &\le& (c_4+c_5) \max\{|\beta_1|^{\delta_b m}, |\alpha_1|^{\delta_a n}\}+ c_3^{1/2} |\alpha_1|^{n-c_0\log^2 n}. \nonumber
\end{eqnarray}

We distinguish two cases depending on the maximum.

\medskip
{\bf Case I. $|\alpha_1|^{\delta_a n} \ge |\beta_1|^{\delta_b m}$.} The first summand on the right hand side of \eqref{e:kombinalt} is now at most $(c_4+c_5)|\alpha_1|^{\delta_a n}$. Dividing by $|\alpha_1|^n$ and taking into account that $\delta_a<1$ we get
$$
\left| \frac{A_{n,m}}{\alpha_1^n}\right| < (c_4+c_5)|\alpha_1|^{(\delta_a-1) n} + c_3^{1/2} |\alpha_1|^{-c_0\log^2 n} < 2c_3^{1/2} |\alpha_1|^{-c_0\log^2 n} < |A_1|/2,
$$
whenever $n$ is large enough.
This implies
\begin{equation}\label{e:kombinalt1}
|A_1| |\alpha_1|^n/2 < |B_1 \beta_1^m + B_2 \beta_2^m| < 3|A_1| |\alpha_1|^n/2.
\end{equation}
A direct consequence of this inequality and \eqref{e:komplexkorlat1} is that
$$
|A_1| |\alpha_1|^n/2 < c_7 |\beta_1|^m,
$$
thus
\begin{equation}\nonumber
  n- \delta m < c_8 \quad \mbox{with}\quad c_8 = \frac{\log c_7 - \log (|A_1|/2) }{\log |\alpha_1|} .
\end{equation}

On the other hand he inequalities \eqref{e:komplexkorlat1} and \eqref{e:kombinalt1} imply, 
$$
3|A_1| |\alpha_1|^n/2> |\beta_1|^{m-c_6 \log m},
$$
thus
$$
n- \delta m > - c_6 \delta \log m - c_9 \quad \mbox{with} \quad c_9 = \frac{\log(3|A_1|/2)}{\log|\alpha_1|}.
$$

\medskip
{\bf Case II. $|\alpha_1|^{\delta_a n} < |\beta_1|^{\delta_b m}$.} The first summand of the right hand side of \eqref{e:kombinalt} is now at most $(c_4+c_5)|\beta_1|^{\delta_b m}$. Our assumption ($a_n\ge b_m$) together with the previous estimations implies
$$
|\alpha_1|^n> \frac{1}{2|A_1|}|a_n|\ge \frac{1}{2|A_1|}|b_m| > \frac{1}{2|A_1|}|\beta_1|^{m-2c_6 \log m}.
$$
Dividing \eqref{e:kombinalt} by $|\alpha_1|^n$ we get
\begin{eqnarray} \label{e:kozepso}
 \label{.} \left| A_1 - \frac{B_1 \beta_1^m + B_2 \beta_2^m}{\alpha_1^n}\right|  &<& (c_4+c_5)\frac{|\beta_1|^{\delta_b m}}{|\alpha_1|^n} + c_3^{1/2} |\alpha_1|^{-c_0\log^2 n}  \nonumber\\
   &<& 2(c_4+c_5)|A_1||\beta_1|^{\delta_b m-m+2c_6 \log m} +c_3^{1/2} |\alpha_1|^{-c_0\log^2 n} \\
   &<& |A_1|/2, \nonumber
\end{eqnarray}
provided $m$ (thus $n$) is large enough. From here on we can repeat the argument of Case I and get the same lower and upper estimates for $n-\delta m$.

For later use we summarize them in the inequality
\begin{equation} \label{e:n-m delta}
 -c_{10} \log m < n- \delta m <  c_8,
\end{equation}
where $c_{10}= c_4\delta + c_9$ and $c_8$ are positive constants. Its simpler, but weaker form is,
\begin{equation} \label{e:n-m delta1}
\delta m/2 < n < 2\delta m,
\end{equation}
which is sometimes easier to apply and is valid if $m > \max\{c_8/(2\delta), (2c_{10}/\delta)^2\}$. Plainly, if $m$ does not satisfy the last inequality, but $(n,m)$ is a solution of \eqref{e:kisebb} $n$ is effectively bounded too.

In Case I we proved $|A_{n,m}|< 2c_3^{1/2} |\alpha_1|^{n-c_0\log^2 n}$. Now we are in the position to prove the same bound for Case II too. Indeed, we can estimate the exponent of the first term of \eqref{e:kozepso}
$$
\delta_b m-m+2c_6 \log m< \frac{\delta_b-1}{2\delta}n+ 2c_6 \log(2n/\delta)< -c_{11}n
$$
with an effective, positive constant $c_{11}$, which is our claim.

So far we proved that if $(n,m)\in \Z^2$ is a solution of \eqref{e:kisebb} such that $n$ is large enough then
\begin{equation}\label{e:summa}
  |A_{n,m}| <2c_3^{1/2} |\alpha_1|^{n-c_0\log^2 n}
\end{equation}
and \eqref{e:n-m delta} holds.

\medskip
In the next step we prove a lower bound for $|A_{n,m}|$. To achieve that write $\beta_1= rz$ with $r\in \R_{>0}$ and $z\in \C, |z|=1$, and we see that $\beta_2=r/z$ and $B_1 \beta_1^m + B_2 \beta_2^m = r^m(B_1 z^m + B_2 z^{-m})$. Remember that $r=|\beta_1|=|\alpha_1|^{\delta}$. Substituting this into \eqref{e:kombinalt}, dividing  by $|B_1r^m|$ and using \eqref{e:n-m delta} and \eqref{e:summa} we get
\begin{eqnarray}
  \left| z^m - \frac{A_1\alpha_1^n}{B_1 r^m} + \frac{B_2}{B_1}z^{-m}\right| &\le & \frac{2c_3^{1/2}}{|B_1|}\frac{|\alpha_1|^{n-c_0\log^2 n}}{r^m}\nonumber\\
   &=& \frac{2c_3^{1/2}}{|B_1|}|\alpha_1|^{n-\delta m-c_0\log^2 n} \label{e:fkorlat}\\
   &<& c_{12}|\alpha_1|^{-c_0\log^2 n}, \nonumber
\end{eqnarray}
with $c_{12} = 2|\alpha_1|^{c_8}c_3^{1/2}/|B_1|$ .

The next argumentation was used for the tribonacci numbers by Bravo et al. \cite{BGL, BGLTK} (especially the proof of Lemma 1 in \cite{BGLTK}).
The zeroes of the quadratic polynomial $X^2 - \frac{A_1\alpha_1^n}{B_1 r^m} X + \frac{B_2}{B_1}$ are
$$
2\lambda_1 =  \frac{A_1\alpha_1^n}{B_1r^m} + \sqrt{\left( \frac{A_1\alpha_1^n}{B_1r^m} \right)^2 -4 \frac{B_2}{B_1}}, \quad 2\lambda_2 =  \frac{A_1\alpha_1^n}{B_1r^m} - \sqrt{\left( \frac{A_1\alpha_1^n}{B_1r^m} \right)^2 -4 \frac{B_2}{B_1}},
$$
thus
$$
|z^{m} - \lambda_1||z^{m} - \lambda_2| \le c_{12}|\alpha_1|^{-c_0\log^2 n}.
$$
We may assume without loss of generality $|z^{m} - \lambda_1|\le |z^{m} - \lambda_2|$, hence
$$
|z^{m} - \lambda_1|\le \sqrt{c_{12}} |\alpha_1|^{-c_0(\log^2 n)/2},
$$
which implies
$$
|\lambda_1|\ge |z|^m - \sqrt{c_{12}} |\alpha_1|^{-c_0(\log^2 n)/2} \ge 1/2,
$$
provided $n$ is large enough. Thus
\begin{equation}\label{e:felső}
  |z^m/\lambda_1 - 1| \le 2 \sqrt{c_{12}} |\alpha_1|^{-c_0(\log^2 n)/2}.
\end{equation}

Set $\Gamma= z^m/\lambda_1 -1$. If $\Gamma=0$ then $A_1 \alpha^n = B_1\beta_1^m + B_2\beta_2^m$. For an algebraic integer $\eta$ denote $P(\eta)$ the largest prime lying below the prime ideal divisors of $(\eta)$. On one hand $P(A_1 \alpha^n)$ is uniformly bounded, but as $\beta_1/\beta_2$ is not a root of unity $P(B_1\beta_1^m + B_2\beta_2^m)$ can be arbitrarily large by Peth\H o \cite{P0}. Thus $\Gamma\not=0$ and we can apply Lemma \ref{Matveev} with he actual parameters:
$$
\K= \Q(z,\lambda_1), t=2, \eta_1= z, \eta_2 = \lambda_1, b_1=m, b_2=1.
$$

As $z=r/\beta_2 = \sqrt{\beta_1/\beta_2}$, the field $\Q(z)$ is included in a quadratic extension of $Q(\beta_1,\ldots,\beta_l)$. Similarly $\Q(\lambda_1)$ is a subfield of a quadratic extension of $\Q(\alpha_1,\ldots,\alpha_k,\beta_1,\ldots,\beta_l)$, thus $d_{\K}\le 4 k! l!$.

In the preceding argumentation we assumed $n$ large enough several times. Thus, by \eqref{e:n-m delta1} $m$ is large enough too, and we may assume $B=m$. Plainly $h(\eta_1)=h(z)\le 2h(\beta_1)$. The estimation of $h(\eta_2)$ is much more complicated and seems to be possible only if $\delta = \log |\beta_1|/ \log |\alpha_1|$ is rational, which is assumed in the theorem. Indeed, in this case $u = n-m\delta$ is rational too, and as $\frac{\alpha_1^n}{r^m} = \pm |\alpha_1|^u$, we have
$$
h\left( \frac{A_1}{B_1}\frac{\alpha_1^n}{r^m}\right) \le |u| h(\alpha_1)+ h\left( \frac{A_1}{B_1}\right)
$$
by Lemma \ref{l:Weil height} (3). Applying several times Lemma \ref{l:Weil height}  we finally obtain
$$
h(\lambda_1) \le 2|u|h(\alpha_1) + c_{13}
$$
with a computable constant $c_{13}$. Finally taking into account \eqref{e:n-m delta} we get
$$
h(\eta_2)= h(\lambda_1) \le c_{14} \log m + c_{13}
$$
with $c_{14}=  2 h(\alpha_1) c_{10}$. Setting $A_1= 8 k! l! h(\beta_1), A_2= 4 k! l! c_{14} \log m + 4 k! l! c_{13}$ the application of Lemma \ref{Matveev} yields
$$
\log |\Gamma| > - c_{15} \log^2 m - c_{16},
$$
where
$$
c_{15} = 3^{6.5} \cdot 2^{11} \cdot 30^6 (k! l!)^4 \log(k! l!) h(\beta_1) c_{14}
$$
and $c_{16} = c_{15} c_{13}/c_{14}$.

Comparing the lower bound for $\log |\Gamma|$ with the upper bound for $|\Gamma|$, i.e. with \eqref{e:felső} we obtain
$$
c_0 (\log |\alpha_1|) (\log^2 n)/2 -\log(2\sqrt{c_{11}}) < c_{15} \log^2 m + c_{16}.
$$
Using \eqref{e:n-m delta} we get
$$
c_0 < \frac{2c_{14}}{\log |\alpha_1|}\left(1+ \frac{\log(2/\delta)}{\log n}\right)^2 + \frac{2(c_{15}+\log(2\sqrt{c_{11}})}{\log |\alpha_1|} \frac{1}{\log^2 n}
$$
after some simple manipulation. If
$$
n > \max\left\{\frac{2}{\delta}, \exp \sqrt{\frac{c_{16}+\log(2\sqrt{c_{12}})}{c_{15}}} \right\}
$$
then
$$
c_0 < \frac{9c_{15}}{\log |\alpha_1|},
$$
but this contradicts our assumption.

\section{Proof of Theorem \ref{t:masodik}}

Assume that $A(X),B(X)\in \Z[X]$ satisfy the assumptions of Theorem \ref{t:main}. Let $P_1(X),P_2(X)\in \Z[X]$ be such that $\overline{|P_1|}< \overline{|A|}$ and $\overline{|P_2|}< \overline{|B|}$ then $A_1X),B_1(X)$ with $A_1(X)=A(X)P_1(X), B_1(X)=B(X)P_2(X)$ satisfy obviously the assumptions of Theorem \ref{t:main} too. Thus it is enough to characterize the irreducible polynomials among them.

Let $A(X),B(X)\in \Z[X]$ be irreducible, and denote $\alpha$ the real root of $A(X)$ with maximal size, while $\beta$ the non-real root of $B(X)$ with maximal size. In Theorem \ref{t:main} a crucial (and unnatural) assumption is that $|\alpha|$ and $|\beta|$ are multiplicatively dependent. This means that there are coprime, non-zero integers $p,q$ such that $\pm \alpha^p = (\beta \bar{\beta})^{q/2}$ or, equivalently
\begin{equation}\label{e:multdep}
  \alpha^{2p} = (\beta \bar{\beta})^q.
\end{equation}
As $\alpha$ is the dominant real root of $A(X)$ and $\beta$ one of the dominant complex roots of $B(X)$ both have to lie outside the unit circle and the last equation is only possible if $p$ and $q$ have same sign, but then we may assume that both are positive. Denote $\K$ the normal closure of $\Q(\alpha,\beta,\bar{\beta})$ and $G$ its Galois group. If $\sigma\in G$ then $\sigma(\alpha)$ is a root of $A(X)$, thus $|\sigma(\alpha)| < |\alpha|$. Because of the same reason $|\sigma(\beta)| < |\beta|$ for all $\sigma \in G$ except when $\sigma(\beta) = \bar{\beta}$.

Let $\sigma \in G$ not the identity be such that $\sigma(\alpha) =\alpha$. If $\sigma(\beta) \not= \bar{\beta}$ then
$$
(\beta \bar{\beta})^q =\alpha^{2p} = \sigma(\alpha^{2p}) = \left(\sigma(\beta) \sigma(\bar{\beta})\right)^q,
$$
but this is impossible because $|\sigma(\beta)| < |\beta|$ and $|\sigma(\bar{\beta})| \le |\bar{\beta}|$. Thus either $\alpha \in \Z$ and $\deg G = 2$ or $B(X)$ has apart $\bar{\beta}$ a third root, say $\gamma$.

If $\alpha\in \Z$ then the only algebraic conjugate of $\beta$ is $\bar{\beta}$, thus $\beta\bar{\beta}\in \Z$, which has to be a $p$-th power, say $b^p$. Hence $\alpha^2=b^q$ and $\alpha=\pm b^{q/2}$. If $b$ is a square then $q$ is arbitrary, otherwise it must be even. Hence the defining polynomial of $\beta$ has the shape $X^2+aX+b^p$ such that $a^2-4b^p<0$, further the defining polynomial of $\alpha$ is $X\pm b^{q/2}$. This proves that if $A$ has an integer root then (i) is not only sufficient, but also necessary that the pair $A(X),B(X)$ fulfill the requirements of Theorem \ref{t:main}.

We now turn to the second case, but now we can prove only sufficiency. Denote $\alpha$ the real root and $\gamma, \overline{\gamma}$ the pair of conjugate complex roots of $Q_2(X)$. We have $\alpha \gamma \overline{\gamma} =1$, thus $|\gamma|<1$ because by assumption $|\alpha|>1$. Further the last equation means  $|\alpha| |\gamma|^2 = 1$ too, i.e $|\alpha|$ and $|\gamma|$ are multiplicatively dependent. The roots of $X^3Q_2(1/X)$ are $1/\alpha, 1/\gamma, 1/\overline{\gamma}$. With $\beta = 1/\gamma$ we have that $\alpha$ and $\beta$ are multiplicatively dependent. 

It remains to prove that neither $\alpha/\beta$ nor $\beta/\overline{\beta}$ is a root of unity. In contrary assume that $\alpha/\beta= \omega$ is a root of unity. Taking conjugate and multiplying the resulting equation with the first one we get
$$
\frac{\alpha}{\beta} \frac{\alpha}{\overline{\beta} } = \omega \overline{\omega} =1.
$$
We also have $\alpha/(\beta\overline{\beta})=1 $, which together with the last equation implies $\alpha^3 = 1$, i.e. $\alpha=1$, which is absurd.

Assume finally that $\beta/\overline{\beta} = \omega$ is a root of unity. Then, especially $\omega \in \K= \Q(\beta, \overline{\beta})$. We also have $\alpha \in \K$ because $\alpha + \beta + \overline{\beta} = a$, thus the field $\K$ is normal. There is an automorphism $\sigma$ of $\K$, with $\sigma(\beta)=\alpha$. It maps $\overline{\beta})$ to itself or to $\beta$ and $\omega$ to one of its conjugates, which is also root of unity. In both cases $\alpha/\beta$ is a root of unity, which is by the last paragraph absurd. Hence the assumptions of Theorem \ref{t:main} fulfill. $\Box$

\section{Proof of Corllary \ref{c:Bravo}}

The characteristic polynomial of $(f_n)$ is $F(X)=X^3-aX^2-bX-1$, while that of $(f_{-n})$ is $F_{-}(X) = X^3+bX^2+aX-1 = -X^3 F(\frac{1}{X})$, hence the roots of $F(X)$ and $F_{-}(X)$ are reciprocal. We have further $F(X) = -((-X)^3 +a(-X)^2-b(-X)+1)$. The polynomial in bracket has the same shape as $Q_2(X)$, but with $-b$ instead of $b$. The assumption $-27-18ab+a^2b^2-4a^3+4b^3<0$ means that its discriminant is negative, i.e. $F(X)$ has one real and a pair of conjugate complex roots. Denote by $\alpha$ its real root. Because $F(0)=-1$ and $F(x) \to \infty$ with $x$, $\alpha$ has to be positive. Interchanging the positive and negative directions we may assume $\alpha>1$. Denote the other roots of $F(X)$ by $1/\beta, 1/\overline{\beta}$. They are in absolute value less than one. The roots of $F_{-}(X)$ are $\beta,\overline{\beta}, 1/\alpha$. With these notations we have 
\begin{equation} \label{e:Binet}
f_n = a \alpha^n + b \beta^n + \overline{b} \overline{\beta}^n, \; n\in \Z
\end{equation}
with some non-zero $a,b\in \C$. To prove the corollary we have to distinguish three cases.

\medskip
{\bf Case 1, $n,m\ge 0$.} A simple, elementary argument, which is left to the reader, is enough to settle it.

\medskip
{\bf Case 2, $n,m< 0$.} Now $f_{-n}$ is still exponentially growing, but also oscillating. To prove the corollary we have to invoke the theory of linear forms in logarithms of algebraic numbers. More precisely Theorem 3.5 of \cite{ShT} with the choice $x_n=y_n = b \beta^n + \overline{b} \overline{\beta}^n, A=B=1$ implies the statement. 

\medskip
{\bf Case 3, $n\cdot m< 0$.} We may assume $n>0, m<0$ without loss of generality. Setting $m=-k, k>0$ our equation becomes
$$
f_n = f_{-k}.
$$
The characteristic polynomial of $(f_n)$ is $F(X)$ and that of $(f_{-n})$ is $F_{-}(X)$. By the remarks at the beginning of the proof setting $Q_2(X)=-F(-X), P_1(X)=P_2(X)=1$ the assumptions of Theorem \ref{t:masodik} hold, i.e. $(f_n)$ and $(f_{-n})$ are such sequences for which Corollary \ref{c:main} is applicable, thus our equation has only finitely many effectively computable solutions in this case too..

\end{document}